\documentclass{amsproc}%
\usepackage{graphicx}
\usepackage{amscd}
\usepackage{amsmath}
\usepackage{amsfonts}
\usepackage{amssymb}%
\theoremstyle{plain}

\numberwithin{equation}{section}

\begin{document}
\title[The centralizer of an element in an endomorphism ring ]{The centralizer of an element in an endomorphism ring }
\author{Jen\H{o} Szigeti}
\address{Institute of Mathematics, University of Miskolc, Miskolc, Hungary 3515}
\email{jeno.szigeti@uni-miskolc.hu}
\author{Leon van Wyk}
\address{Department of Mathematical Sciences, Stellenbosch University\\
P/Bag X1, Matieland 7602, South Africa }
\email{lvw@sun.ac.za}
\thanks{\noindent The first author was supported by OTKA of Hungary No. K61007. The
second author was supported by the National Research Foundation of South
Africa under Grant No. UID 61857. Any opinion, findings and conclusions or
recommendations expressed in this material are those of the authors and
therefore the National Research Foundation does not accept any liability in
regard thereto.}

\begin{abstract}
We prove that the centralizer $C_{\varphi}\subseteq$ Hom$_{R}(M,M)$ of a
nilpotent endomorphism $\varphi$\ of a finitely generated semisimple left
$R$-module $_{R}M$ (over an arbitrary ring $R$) is the homomorphic image of
the opposite of a certain $Z(R)$-subalgebra of the full $m\times m$ matrix
algebra $M_{m\times m}(R[t])$, where $m$ is the dimension (composition length)
of $\ker(\varphi)$. If $R$ is a finite dimensional division ring over its
central subfield $Z(R)$ and $\varphi$ is nilpotent, then we give an upper
bound for the $Z(R)$-dimension of $C_{\varphi}$. If $R$ is a local ring,
$\varphi$ is nilpotent and $\sigma\in$Hom$_{R}(M,M)$ is arbitrary, then we
provide a complete description of the containment $C_{\varphi}\subseteq
C_{\sigma}$ in terms of an appropriate $R$-generating set of $_{R}M$. For an
arbitrary (not necessarily nilpotent) linear map $\varphi\in$Hom$_{K}(V,V)$ of
a finite dimensional vector space $V$\ over an algebraically closed field $K$
we prove that $C_{\varphi}$ is the homomorphic image of a direct product of
$p$ factors such that for each $1\leq i\leq p$ the $i$-th factor is a
$K$-subalgebra of $M_{m_{i}\times m_{i}}(K[t])$ with $m_{i}=\dim(\ker
(\varphi-\lambda_{i}1_{V}))$ and $\{\lambda_{1},\lambda_{2},\ldots,\lambda
_{p}\}$ is the set of all eigenvalues of $\varphi$. As a consequence, we
obtain that $C_{\varphi}$\ satisfies all polynomial identities of $M_{m\times
m}(K[t])$, where $m$ is the maximum of the $m_{i}$'s.

\end{abstract}
\maketitle

\noindent1. INTRODUCTION

\bigskip

Our work was motivated by one of the classical subjects of advanced linear
algebra. A detailed study of commuting matrices can be found in many of the
text books on linear algebra ([5,6]). Commuting pairs and $k$-tuples of
$n\times n$\ matrices have been continuously in the focus of research (see,
for example, [2,3,4]). If we replace $n\times n$ matrices by endomorphisms of
an $n$-generated module, we get a more general situation. The aim of the
present paper is to investigate the size and the PI properties of the
centralizer $C_{\varphi}$\ of an element $\varphi$ in the endomorphism ring
Hom$_{R}(M,M)$.

In the case of a finitely generated semisimple left $R$-module $_{R}M$ (over
an arbitrary ring $R$) we obtain results about the centralizer of a nilpotent
endomorphism $\varphi$. First we prove that $C_{\varphi}$ is the homomorphic
image of the opposite of a certain $Z(R)$-subalgebra of the full $m\times m$
matrix algebra $M_{m\times m}(R[t])$\ over the polynomial ring $R[t]$, where
$m$ is the dimension (composition length) of $\ker(\varphi)$. As a
consequence, we obtain that $C_{\varphi}$\ satisfies all polynomial identities
of $M_{m\times m}^{\text{op}}(R[t])$; in particular, if $R$ is commutative,
then the standard identity $S_{2m}=0$ of degree $2m$ holds in $C_{\varphi}$.
The authors do not know similar results in the literature.

Then we exhibit the centralizer $C_{\varphi}$\ of a nilpotent $\varphi$ as a
homomorphic image of the intersection of some $md$-generated $R$-submodule of
$M_{m\times m}(R[t])$ and of the above mentioned $Z(R)$-subalgebra, where $d$
is the dimension (composition length) of $_{R}M$. If $R$ is a finite
dimensional division ring over its central subfield $Z(R)$, then we can use
this homomorphic image representation to derive a \textquotedblleft
sharp\textquotedblright\ upper bound for the $Z(R)$-dimension of $C_{\varphi}%
$. For an arbitrary (not necessarily nilpotent) linear map of a finite
dimensional vector space, or equivalently for a matrix $A\in M_{n\times n}(K)$
over an algebraically closed field $K$, the dimension of $C_{A}$ is well known
(see [5,6]). A multiplicative ($K$ vector space) base of $C_{A}$ was
constructed in [4].

If $\varphi:M\longrightarrow M$ is a so called indecomposable nilpotent
$R$-endomorphism, then we give a complete description of $C_{\varphi}$ in
terms of an appropriate $R$-generating set of $_{R}M$. In particular, if $R$
is commutative, then we prove that $\psi\in C_{\varphi}$ holds if and only if
$\psi$ is a polynomial of $\varphi$. A nilpotent linear map (over an
algebraically closed field) is indecomposable if and only if its
characteristic polynomial coincides with the minimal polynomial. The following
is a classical result about the centralizer (see [5] part VII, section 39). If
$K$\ is an algebraically closed field and the characteristic polynomial of a
(not necessarily nilpotent) $A\in M_{n\times n}(K)$ coincides with the minimal
polynomial, then%
\[
C_{A}=\{f(A)\mid f(t)\in K[t]\}.
\]
For an indecomposable nilpotent $\varphi$\ our description of $C_{\varphi}$ is
a generalization of the above result. We note that the above result on $C_{A}%
$\ is similar to Bergman's Theorem ([1]) about the centralizer of a
non-constant polynomial in the free associative algebra.

The containment relation $C_{\varphi}\subseteq C_{\sigma}$ is also considered.
If $R$ is a local ring, $\varphi$ is nilpotent and $\sigma\in$ Hom$_{R}(M,M)$
is arbitrary, then we provide a complete description of the situation
$C_{\varphi}\subseteq C_{\sigma}$ in terms of an appropriate $R$-generating
set of $_{R}M$. For two (not necessarily nilpotent) matrices $A,B\in
M_{n\times n}(K)$, over an algebraically closed field $K$, the containment
$C_{A}\subseteq C_{B}$ holds if and only if $B=f(A)$ for some $f(t)\in K[t]$
(see [5] part VII, section 39). For a nilpotent $\varphi$\ our description of
$C_{\varphi}\subseteq C_{\sigma}$ is a generalization of the above result.

In the case of a finite dimensional vector space $V$ over an algebraically
closed field $K$ we obtain results about the centralizer of an arbitrary (not
necessarily nilpotent) linear map $\varphi\in$ Hom$_{K}(V,V)$. We prove that
$C_{\varphi}$ is the homomorphic image of a direct product of $p$ factors such
that for each $1\leq i\leq p$ the $i$-th factor is a $K$-subalgebra of
$M_{m_{i}\times m_{i}}(K[t])$, where%
\[
m_{i}=\dim(\ker(\varphi-\lambda_{i}1_{V}))
\]
and $\{\lambda_{1},\lambda_{2},\ldots,\lambda_{p}\}$ is the set of all
eigenvalues of $\varphi$. As a consequence, we obtain that $C_{\varphi}%
$\ satisfies all polynomial identities of $M_{m\times m}(K[t])$, where $m$ is
the maximum of the $m_{i}$'s. The following related results can be found in
[4,6]. The centralizer $K$-algebra $C_{A}$\ of a matrix $A\in M_{n\times
n}(K)$ (over an algebraically closed field $K$) is isomorphic to the direct
product of the centralizers $C_{A_{i}}$, where $A_{i}$ denotes the block
diagonal matrix consisting of all Jordan blocks of $A$ having eigenvalue
$\lambda_{i}$ in the diagonal. The number of the diagonal blocks in $A_{i}$ is
$m_{i}=\dim(\ker(A-\lambda_{i}I))$ and the size of $A_{i}$ is $d_{i}\times
d_{i}$, where $d_{i}$ is the multiplicity of the root $\lambda_{i}$\ in the
characteristic polynomial of $A$. For a matrix $X\in M_{d_{i}\times d_{i}}(K)$
we have $X\in C_{A_{i}}$ if and only if $X=[X_{k,l}]$ is an $m_{i}\times
m_{i}$ matrix of blocks, where $X_{k,l}$ is an arbitrary triangularly striped
matrix block of size $s_{i}(k)\times s_{i}(l)$ and $s_{i}(k)$ denotes the size
of the $k$-th elementary Jordan block in $A_{i}$ for $1\leq k\leq m_{i}$. We
note that $\dim(C_{A_{i}})=s_{i}(1)+3s_{i}(2)+\cdots+(2m_{i}-1)s_{i}(m_{i})$
provided that $s_{i}(1)\geq s_{i}(2)\geq\cdots\geq s_{i}(m_{i})$. It seems
that we cannot use the striped block structure of the elements of $C_{A_{i}}$
to derive the above mentioned PI properties of $C_{A}$ ($\cong C_{\varphi}$),
only the much weaker statement that $C_{\varphi}$\ (or $C_{A}$) satisfies all
polynomial identities of $M_{d\times d}(K)$ follows ($d$ is the maximum of the
$d_{i}$'s and $d_{i}=s_{i}(1)+s_{i}(2)+\cdots+s_{i}(m_{i})$). The above
consideration confirms that the use of the polynomial ring $K[t]$ is an
essential ingredient of our treatment.

Since all known results about centralizers are in close connection with the
Jordan normal base, it is not surprising that our development depends on the
existence of the (nilpotent) Jordan normal base of a semisimple module with
respect to a given nilpotent endomorphism (guaranteed by one of the main
theorems of [7]). Using an endomorphism of $_{R}M$, a natural $R[t]$-module
structure on $M$ can be defined. Almost all of our proofs are based on
calculations using this induced module structure. In the proof of the last
theorem we use the Fitting Lemma.

\bigskip

\noindent2. THE NILPOTENT\ JORDAN NORMAL\ BASE

\bigskip

\noindent Throughout the paper a ring $R$\ means a (not necessarily
commutative) ring with identity, $Z(R)$ and $J(R)$ denote the centre and the
Jacobson radical of $R$, respectively. Let $\varphi:M\longrightarrow M$ be an
$R$-endomorphism of the (unitary) left $R$-module $_{R}M$. A subset%
\[
\{x_{\gamma,i}\mid\gamma\in\Gamma,1\leq i\leq k_{\gamma}\}\subseteq M
\]
is called a \textit{nilpotent Jordan normal base} of $_{R}M$\ with respect to
$\varphi$ if each $R$-submodule $Rx_{\gamma,i}\leq M$ is simple,%
\[
\underset{\gamma\in\Gamma,1\leq i\leq k_{\gamma}}{\oplus}Rx_{\gamma,i}=M
\]
is a direct sum, $\varphi(x_{\gamma,i})=x_{\gamma,i+1}$ $,$ $\varphi
(x_{\gamma,k_{\gamma}})=x_{\gamma,k_{\gamma}+1}=0$ for all\textit{ }$\gamma
\in\Gamma$, $1\leq i\leq k_{\gamma}$, and the set $\{k_{\gamma}\mid\gamma
\in\Gamma\}$ of integers is bounded. For $i\geq k_{\gamma}+1$ we assume that
$x_{\gamma,i}=0$ holds in $M$. Now $\Gamma$\ is called the set of (Jordan-)
blocks and the size of the block $\gamma\in\Gamma$ is the integer $k_{\gamma
}\geq1$. Obviously, the existence of a nilpotent Jordan normal base implies
that $_{R}M$ is semisimple and $\varphi$ is nilpotent with $\varphi^{n}%
=0\neq\varphi^{n-1}$, where%
\[
n=\max\{k_{\gamma}\mid\gamma\in\Gamma\}
\]
is the index of nilpotency. If $_{R}M$ is finitely generated, then%
\[
\underset{\gamma\in\Gamma}{\sum}k_{\gamma}=\dim_{R}(M)
\]
is the dimension of $_{R}M$\ (equivalently: the composition length of $_{R}M$
or the height of the submodule lattice of $_{R}M$). Clearly,%
\[
\varphi(M)=\varphi\left(  \underset{\gamma\in\Gamma,1\leq i\leq k_{\gamma}%
}{\sum}Rx_{\gamma,i}\right)  =\underset{\gamma\in\Gamma,1\leq i\leq k_{\gamma
}}{\sum}R\varphi(x_{\gamma,i})=\underset{\gamma\in\Gamma,1\leq i\leq
k_{\gamma}-1}{\sum}Rx_{\gamma,i+1}%
\]
implies that%
\[
\text{im}(\varphi)=\underset{\gamma\in\Gamma,1\leq i\leq k_{\gamma}-1}{\oplus
}Rx_{\gamma,i+1}=\underset{\gamma\in\Gamma^{\prime},2\leq i^{\prime}\leq
k_{\gamma}}{\oplus}Rx_{\gamma,i^{\prime}},
\]
where%
\[
\Gamma^{\prime}=\{\gamma\in\Gamma\mid k_{\gamma}\geq2\}\text{ and }%
\Gamma\setminus\Gamma^{\prime}=\{\gamma\in\Gamma\mid k_{\gamma}=1\}.
\]
Any element $u\in M$ can be written as%
\[
u=\underset{\gamma\in\Gamma,1\leq i\leq k_{\gamma}}{\sum}a_{\gamma,i}%
x_{\gamma,i},
\]
where $\{(\gamma,i)\mid\gamma\in\Gamma,1\leq i\leq k_{\gamma},$ and
$a_{\gamma,i}\neq0\}$ is finite and all summands $a_{\gamma,i}x_{\gamma,i}$
are uniquely determined by $u$. Since%
\[
\varphi(u)=\underset{\gamma\in\Gamma,1\leq i\leq k_{\gamma}}{\sum}a_{\gamma
,i}\varphi(x_{\gamma,i})=\underset{\gamma\in\Gamma,1\leq i\leq k_{\gamma}%
}{\sum}a_{\gamma,i}x_{\gamma,i+1}=0
\]
is equivalent to the condition that $a_{\gamma,i}x_{\gamma,i+1}=0$ for all
$\gamma\in\Gamma,1\leq i\leq k_{\gamma}-1$, we obtain that%
\[
\varphi(u)=0\Longleftrightarrow u=\underset{\gamma\in\Gamma}{\sum}%
a_{\gamma,k_{\gamma}}x_{\gamma,k_{\gamma}}.
\]
Indeed, $a_{\gamma,i}x_{\gamma,i}\neq0$ ($1\leq i\leq k_{\gamma}-1$) would
imply that $ba_{\gamma,i}x_{\gamma,i}=x_{\gamma,i}$ for some $b\in R$ (note
that $Rx_{\gamma,i}\leq M$ is simple), whence%
\[
x_{\gamma,i+1}=\varphi(x_{\gamma,i})=ba_{\gamma,i}\varphi(x_{\gamma
,i})=ba_{\gamma,i}x_{\gamma,i+1}=0
\]
can be derived, a contradiction. It follows that%
\[
\ker(\varphi)=\underset{\gamma\in\Gamma}{\oplus}Rx_{\gamma,k_{\gamma}}%
\]
and $\dim_{R}(\ker(\varphi))=\left\vert \Gamma\right\vert $ in case of a
finite $\Gamma$. The following is one of the main results in [7].

\bigskip

\noindent\textbf{2.1.Theorem.}\textit{ Let }$\varphi:M\longrightarrow
M$\textit{ be an }$R$\textit{-endomorphism of the left }$R$\textit{-module
}$_{R}M$\textit{. Then the following are equivalent.}

\begin{enumerate}
\item $_{R}M$\textit{ is a semisimple left }$R$\textit{-module and }$\varphi
$\textit{\ is nilpotent.}

\item \textit{There exists a nilpotent Jordan normal base of }$_{R}%
M$\textit{\ with respect to }$\varphi$\textit{.}
\end{enumerate}

\bigskip

\noindent\textbf{2.2.Proposition.}\textit{ Let }$\varphi:M\longrightarrow
M$\textit{ be a nilpotent }$R$\textit{-endomorphism of the finitely generated
semisimple left }$R$\textit{-module }$_{R}M$\textit{. If }$\{x_{\gamma,i}%
\mid\gamma\in\Gamma,1\leq i\leq k_{\gamma}\}$\textit{ and }$\{y_{\delta,j}%
\mid\delta\in\Delta,1\leq j\leq l_{\delta}\}$\textit{\ are nilpotent Jordan
normal bases of }$_{R}M$\textit{\ with respect to }$\varphi$\textit{, then
there exists a bijection }$\pi:\Gamma\longrightarrow\Delta$\textit{ such that
}$k_{\gamma}=l_{\pi(\gamma)}$\textit{ for all }$\gamma\in\Gamma$\textit{. Thus
the sizes of the blocks of a nilpotent Jordan normal base are unique up to a
permutation of the blocks.}

\bigskip

\noindent\textbf{Proof.} We apply induction on the index of the nilpotency of
$\varphi$. If $\varphi=0$, then we have $k_{\gamma}=l_{\delta}=1$ for all
$\gamma\in\Gamma$, $\delta\in\Delta$, and%
\[
\underset{\gamma\in\Gamma}{\oplus}Rx_{\gamma,1}=\underset{\delta\in\Delta
}{\oplus}Ry_{\delta,1}=M
\]
implies the existence of a bijection $\pi:\Gamma\longrightarrow\Delta$
(Krull-Schmidt, Kuros-Ore). Assume that our statement holds for any
$R$-endomorphism $\phi:N\longrightarrow N$ with $_{R}N$ being a finitely
generated semisimple left $R$-module and $\phi^{n-2}\neq0=\phi^{n-1}$.
Consider the situation described in the proposition with $\varphi^{n-1}%
\neq0=\varphi^{n}$, then%
\[
\text{im}(\varphi)=\underset{\gamma\in\Gamma^{\prime},2\leq i^{\prime}\leq
k_{\gamma}}{\oplus}Rx_{\gamma,i^{\prime}}%
\]
ensures that%
\[
\{x_{\gamma,i^{\prime}}\mid\gamma\in\Gamma^{\prime},2\leq i^{\prime}\leq
k_{\gamma}\}
\]
is a nilpotent Jordan normal base of the left $R$-submodule im$(\varphi)\leq
M$ of $_{R}M$\ with respect to the restricted $R$-endomorphism $\varphi:$
im$(\varphi)\longrightarrow$ im$(\varphi)$. The same holds for%
\[
\{y_{\delta,j^{\prime}}\mid\delta\in\Delta^{\prime},2\leq j^{\prime}\leq
l_{\delta}\}.
\]
Since we have $\phi^{n-2}\neq0=\phi^{n-1}$ for $\phi=\varphi\upharpoonright$
im$(\varphi)$, our assumption ensures the existence of a bijection $\pi
:\Gamma^{\prime}\longrightarrow\Delta^{\prime}$ such that $k_{\gamma}%
-1=l_{\pi(\gamma)}-1$ for all $\gamma\in\Gamma^{\prime}$. In view of%
\[
\ker(\varphi)=\underset{\gamma\in\Gamma}{\oplus}Rx_{\gamma,k_{\gamma}%
}=\underset{\delta\in\Delta}{\oplus}Ry_{\delta,l_{\delta}}%
\]
we obtain that $\left\vert \Gamma\right\vert =\left\vert \Delta\right\vert $
(Krull-Schmidt, Kuros-Ore), whence $\left\vert \Gamma\setminus\Gamma^{\prime
}\right\vert =\left\vert \Delta\setminus\Delta^{\prime}\right\vert $ follows.
Thus we have a bijection $\pi^{\ast}:\Gamma\setminus\Gamma^{\prime
}\longrightarrow\Delta\setminus\Delta^{\prime}$ and the natural map%
\[
\pi\sqcup\pi^{\ast}:\Gamma^{\prime}\cup(\Gamma\setminus\Gamma^{\prime
})\longrightarrow\Delta^{\prime}\cup(\Delta\setminus\Delta^{\prime})
\]
is a bijection with the desired property.$\square$

\bigskip

\noindent We call a nilpotent element $s\in S$ of the ring $S$%
\ \textit{decomposable} if $es=se$ holds for some idempotent element $e\in S$
($e^{2}=e$) with $0\neq e\neq1$. A nilpotent element which is not decomposable
is called \textit{indecomposable}.

\bigskip

\noindent\textbf{2.3.Proposition.}\textit{ Let }$\varphi:M\longrightarrow
M$\textit{ be a nonzero nilpotent }$R$\textit{-endomorphism of the semisimple
left }$R$\textit{-module }$_{R}M$\textit{. Then the following are equivalent.}

\begin{enumerate}
\item \textit{There is a nilpotent Jordan normal base }$\{x_{i}\mid1\leq i\leq
n\}$\textit{\ of }$_{R}M$\textit{\ with respect to }$\varphi$\textit{
consisting of one block (thus }$\left\vert \Gamma\right\vert =1$\textit{\ for
any nilpotent Jordan normal base }$\{x_{\gamma,i}\mid\gamma\in\Gamma,1\leq
i\leq k_{\gamma}\}$\textit{\ of }$_{R}M$\textit{\ with respect to }$\varphi
$\textit{).}

\item $\varphi$\textit{ is an indecomposable nilpotent element of the ring
}Hom$_{R}(M,M)$\textit{.}

\item $_{R}M$\textit{ is finitely generated and }$\varphi^{d-1}\neq0$\textit{,
where }$d=\dim_{R}(M)$\textit{ is the dimension of }$_{R}M$.
\end{enumerate}

\bigskip

\noindent\textbf{Proof.}

\noindent(1)$\Longrightarrow$(3): Clearly,%
\[
\underset{1\leq i\leq n}{\oplus}Rx_{i}=M
\]
implies that we have $d=n$\ for the dimension of $_{R}M$, whence%
\[
\varphi^{d-1}(x_{1})=\varphi^{n-1}(x_{1})=x_{n}\neq0
\]
follows.

\noindent(3)$\Longrightarrow$(1): Let $\{x_{\gamma,i}\mid\gamma\in\Gamma,1\leq
i\leq k_{\gamma}\}$ be a Jordan normal base of $_{R}M$\ with respect to
$\varphi$. Suppose that $\left\vert \Gamma\right\vert \geq2$, then%
\[
n=\max\{k_{\gamma}\mid\gamma\in\Gamma\}\leq d-1,
\]
where $d=\underset{\gamma\in\Gamma}{\sum}k_{\gamma}=\dim_{R}(M)$. Thus
$\varphi^{d-1}=\varphi^{(d-1)-n}\circ\varphi^{n}=0$, a contradiction.

\noindent(1)$\Longrightarrow$(2): Suppose that $\varepsilon\circ
\varphi=\varphi\circ\varepsilon$ holds for some idempotent endomorphism
$\varepsilon\in$Hom$_{R}(M,M)$ with $0\neq\varepsilon\neq1$. Then%
\[
\text{im}(\varepsilon)\oplus\text{im}(1-\varepsilon)=M
\]
for the non-zero (semisimple) $R$-submodules im$(\varepsilon)$ and
im$(1-\varepsilon)$ of $_{R}M$. Now $\varepsilon\circ\varphi=\varphi
\circ\varepsilon$ ensures that $\varphi:$im$(\varepsilon)\longrightarrow$
im$(\varepsilon)$ and $\varphi:$im$(1-\varepsilon)\longrightarrow$
im$(1-\varepsilon)$. Since these restricted $R$-endomorphisms are nilpotent,
we have a nilpotent Jordan normal base of im$(\varepsilon)$ with respect to
$\varphi\upharpoonright$ im$(\varepsilon)$ and a nilpotent Jordan normal base
of im$(1-\varepsilon)$ with respect to $\varphi\upharpoonright$
im$(1-\varepsilon)$. The union of these two bases gives a nilpotent Jordan
normal base of $M$ with respect to $\varphi$ consisting of more than one
block, a contradiction (the direct sum property of the new base is a
consequence of the modularity of the submodule lattice of $_{R}M$).

\noindent(2)$\Longrightarrow$(1): Suppose that $\{x_{\gamma,i}\mid\gamma
\in\Gamma,1\leq i\leq k_{\gamma}\}$ is a nilpotent Jordan normal base of
$_{R}M$\ with respect to $\varphi$ with $\left\vert \Gamma\right\vert \geq2$
and fix an element $\delta\in\Gamma$. Consider the non-zero $\varphi
$-invariant $R$-submodules
\[
N_{\delta}^{\prime}=\underset{1\leq i\leq k_{\delta}}{\oplus}Rx_{\delta
,i}\text{ and }N_{\delta}^{\prime\prime}=\underset{\gamma\in\Gamma
\setminus\{\delta\},1\leq i\leq k_{\gamma}}{\oplus}Rx_{\gamma,i}\text{ },
\]
then $M=N_{\delta}^{\prime}\oplus N_{\delta}^{\prime\prime}$ and define
$\varepsilon_{\delta}:M\longrightarrow M$ as the natural projection of $M$
onto $N_{\delta}^{\prime}$. Then $\varepsilon_{\delta}(u)=u^{\prime}$, where
$u=u^{\prime}+u^{\prime\prime}$ is the unique sum presentation of $u\in M$
with $u^{\prime}\in N_{\delta}^{\prime}$ and $u^{\prime\prime}\in N_{\delta
}^{\prime\prime}$. It is straightforward to see that $\varepsilon_{\delta
}\circ\varepsilon_{\delta}=\varepsilon_{\delta}$, $0\neq\varepsilon_{\delta
}\neq1$\ and $\varepsilon_{\delta}\circ\varphi=\varphi\circ\varepsilon
_{\delta}$ hold.$\square$

\bigskip

\noindent3. THE\ MODULE\ STRUCTURE\ INDUCED\ BY\ AN ENDOMORPHISM

\bigskip

\noindent Let $R[t]$ denote the ring of polynomials of the commuting
indeterminate $t$ with coefficients in $R$. The ideal $(t^{k})=R[t]t^{k}%
=t^{k}R[t]\vartriangleleft R[t]$ generated by $t^{k}$ will be considered in
the sequel. If $\varphi:M\longrightarrow M$ is an arbitrary $R$-endomorphism
of the left $R$-module $_{R}M$, then for $u\in M$ and%
\[
f(t)=a_{1}+a_{2}t+\cdots+a_{n+1}t^{n}\in R[t]
\]
(unusual use of indices!) the left multiplication%
\[
f(t)\ast u=a_{1}u+a_{2}\varphi(u)+\cdots+a_{n+1}\varphi^{n}(u)
\]
defines a natural left $R[t]$-module structure on $M$. This left action of
$R[t]$ on $M$ extends the left action of $R$. The proof of%
\[
g(t)\ast(f(t)\ast u)=(g(t)f(t))\ast u
\]
is straightforward. Note that%
\[
t^{n}\ast u=\varphi^{n}(u)\text{ and }\varphi(f(t)\ast u)=(tf(t))\ast u.
\]
If $\varphi^{k}(u)=0$ for some $1\leq k\leq n$, then%
\[
f(t)\ast u=f^{(k)}(t)\ast u,
\]
where%
\[
f^{(k)}(t)=a_{1}+a_{2}t+\cdots+a_{k}t^{k-1}\in R[t]
\]
is the $k$-cut of $f(t)$. For any $R$-endomorphism $\psi\in$ Hom$_{R}(M,M)$
with $\psi\circ\varphi=\varphi\circ\psi$ we have%
\[
\psi(f(t)\ast u)=f(t)\ast\psi(u)
\]
and hence $\psi:M\longrightarrow M$ is an $R[t]$-endomorphism of the left
$R[t]$-module $_{R[t]}M$. On the other hand, if $\psi:M\longrightarrow M$ is
an $R[t]$-endomorphism of $_{R[t]}M$, then%
\[
\psi(\varphi(u))=\psi(t\ast u)=t\ast\psi(u)=\varphi(\psi(u))
\]
implies that $\psi\circ\varphi=\varphi\circ\psi$. The centralizer
\[
C_{\varphi}=\{\psi\mid\psi\in\text{Hom}_{R}(M,M)\text{ and }\psi\circ
\varphi=\varphi\circ\psi\}
\]
of $\varphi$ is a $Z(R)$-subalgebra of Hom$_{R}(M,M)$ and the argument above
gives that%
\[
C_{\varphi}=\text{Hom}_{R[t]}(M,M).
\]

\noindent For a set $\Gamma\neq\varnothing$, the $\Gamma$-copower
$\underset{\gamma\in\Gamma}{%
{\displaystyle\coprod}
}R[t]$\ of the ring $R[t]$ is an ideal of the $\Gamma$-direct power ring
$\underset{\gamma\in\Gamma}{%
{\displaystyle\prod}
}R[t]$ consisting of all elements $\mathbf{f}=(f_{\gamma}(t))_{\gamma\in
\Gamma}$ with a finite set $\{\gamma\in\Gamma\mid f_{\gamma}(t)\neq0\}$\ of
non-zero coordinates. The power (copower) has a natural $(R[t],R[t])$-bimodule
structure. If $\Gamma$ is finite, then%
\[
(R[t])^{\Gamma}=\underset{\gamma\in\Gamma}{%
{\displaystyle\coprod}
}R[t]=\underset{\gamma\in\Gamma}{%
{\displaystyle\prod}
}R[t].
\]
\noindent If $\{x_{\gamma,i}\mid\gamma\in\Gamma,1\leq i\leq k_{\gamma}\}$ is a
nilpotent Jordan normal base of $_{R}M$\ with respect to a nilpotent
endomorphism $\varphi$, then for an element $\mathbf{f}=(f_{\gamma
}(t))_{\gamma\in\Gamma}$ with%
\[
f_{\gamma}(t)=a_{\gamma,1}+a_{\gamma,2}t+\cdots+a_{\gamma,n_{\gamma}%
+1}t^{n_{\gamma}}%
\]
the formula%
\[
\Phi(\mathbf{f})=\underset{\gamma\in\Gamma,1\leq i\leq k_{\gamma}}{%
{\displaystyle\sum}
}a_{\gamma,i}x_{\gamma,i}=\underset{\gamma\in\Gamma}{%
{\displaystyle\sum}
}f_{\gamma}(t)\ast x_{\gamma,1}%
\]
defines a function%
\[
\Phi:\underset{\gamma\in\Gamma}{%
{\displaystyle\coprod}
}R[t]\longrightarrow M.
\]
In Section 4 we shall make use of%
\[
\Phi(\mathbf{f})=\underset{\gamma\in\Gamma}{%
{\displaystyle\sum}
}f_{\gamma}^{(k_{\gamma})}(t)\ast x_{\gamma,1}=\Phi(\mathbf{f}^{(c)}),
\]
where $f_{\gamma}^{(k_{\gamma})}(t)$\ is the $k_{\gamma}$-cut of $f_{\gamma
}(t)$ and $\mathbf{f}^{(c)}=(f_{\gamma}^{(k_{\gamma})}(t))_{\gamma\in\Gamma}$
(also an element of the copower) is the cut of $\mathbf{f}$ with respect to
the given nilpotent Jordan normal base.

\bigskip

\noindent\textbf{3.1.Theorem.}\textit{ For a nilpotent endomorphism }%
$\varphi\in$Hom$_{R}(M,M)$\textit{\ the function }$\Phi$\textit{\ is a
surjective left }$R[t]$\textit{-homomorphism. We have }$\varphi(\Phi
(\mathbf{f}))=\Phi(t\mathbf{f})$\textit{ for all }$\mathbf{f}\in
\underset{\gamma\in\Gamma}{%
{\displaystyle\coprod}
}R[t]$\textit{ and the kernel}%
\[
\underset{\gamma\in\Gamma}{%
{\displaystyle\coprod}
}J(R)[t]+(t^{k_{\gamma}})\subseteq\ker(\Phi)\vartriangleleft_{l}%
\underset{\gamma\in\Gamma}{%
{\displaystyle\prod}
}R[t]
\]
\textit{is a left ideal of the power (and hence of the copower) ring.}

\bigskip

\noindent\textbf{Proof.} Clearly,%
\[
\underset{\gamma\in\Gamma,1\leq i\leq k_{\gamma}}{\sum}Rx_{\gamma,i}=M
\]
implies that $\Phi$ is surjective. The second part of the defining formula
gives that $\Phi$ is a left $R[t]$-homomorphism:%
\[
\Phi(g(t)\mathbf{f})=\underset{\gamma\in\Gamma}{%
{\displaystyle\sum}
}(g(t)f_{\gamma}(t))\ast x_{\gamma,1}=\underset{\gamma\in\Gamma}{%
{\displaystyle\sum}
}g(t)\ast(f_{\gamma}(t)\ast x_{\gamma,1})=g(t)\ast\Phi(\mathbf{f}),
\]
where $g(t)\in R[t]$. We also have%
\[
\varphi(\Phi(\mathbf{f}))=\underset{\gamma\in\Gamma}{%
{\displaystyle\sum}
}\varphi(f_{\gamma}(t)\ast x_{\gamma,1})=\underset{\gamma\in\Gamma}{%
{\displaystyle\sum}
}(tf_{\gamma}(t))\ast x_{\gamma,1}=\Phi(t\mathbf{f}).
\]
If $\mathbf{f}\in\underset{\gamma\in\Gamma}{%
{\displaystyle\coprod}
}J(R)[t]+(t^{k_{\gamma}})$, then%
\[
f_{\gamma}(t)=(a_{\gamma,1}+a_{\gamma,2}t+\cdots+a_{\gamma,k_{\gamma}%
}t^{k_{\gamma}-1})+(a_{\gamma,k_{\gamma}+1}t^{k_{\gamma}}+\cdots
+a_{\gamma,n_{\gamma}+1}t^{n_{\gamma}})
\]
with $a_{\gamma,i}\in J(R)$, $1\leq i\leq k_{\gamma}$. Since $Rx_{\gamma,i}$
is simple, we have $J(R)x_{\gamma,i}=\{0\}$. Thus $\varphi^{k_{\gamma}%
}(x_{\gamma,1})=0$ implies that $f_{\gamma}(t)\ast x_{\gamma,1}=0$, whence
$\Phi(\mathbf{f})=0$ follows. Take an element $\mathbf{g}=(g_{\gamma
}(t))_{\gamma\in\Gamma}$ of the direct power and suppose that $\Phi
(\mathbf{f})=0$ in $M$. Then%
\[
\underset{\gamma\in\Gamma,1\leq i\leq k_{\gamma}}{\oplus}Rx_{\gamma,i}=M
\]
implies that $a_{\gamma,i}x_{\gamma,i}=0$ for all $\gamma\in\Gamma$ and $1\leq
i\leq k_{\gamma}$. Thus $f_{\gamma}(t)\ast x_{\gamma,1}=0$ for all $\gamma
\in\Gamma$. It follows that%
\[
\Phi(\mathbf{gf})=\underset{\gamma\in\Gamma}{%
{\displaystyle\sum}
}(g_{\gamma}(t)f_{\gamma}(t))\ast x_{\gamma,1}=\underset{\gamma\in\Gamma}{%
{\displaystyle\sum}
}g_{\gamma}(t)\ast(f_{\gamma}(t)\ast x_{\gamma,1})=0,
\]
whence $\mathbf{gf}\in\ker(\Phi)$ can be deduced.$\square$

\bigskip

\noindent If $R$ is a local ring ($R/J(R)$ is a division ring) and
$a_{\gamma,i}x_{\gamma,i}=0$ for some $1\leq i\leq k_{\gamma}$, then
$a_{\gamma,i}\in J(R)$. Thus $\Phi(\mathbf{f})=0$ implies that%
\[
f_{\gamma}(t)=(a_{\gamma,1}+a_{\gamma,2}t+\cdots+a_{\gamma,k_{\gamma}%
}t^{k_{\gamma}-1})+(a_{\gamma,k_{\gamma}+1}t^{k_{\gamma}}+\cdots
+a_{\gamma,n_{\gamma}+1}t^{n_{\gamma}})\in J(R)[t]+(t^{k_{\gamma}}).
\]
It follows that for local rings we have%
\[
\ker(\Phi)=\underset{\gamma\in\Gamma}{%
{\displaystyle\coprod}
}J(R)[t]+(t^{k_{\gamma}}).
\]

\bigskip

\noindent4. THE\ CENTRALIZER\ OF\ A\ NILPOTENT ENDOMORPHISM

\bigskip

\noindent Let $\{x_{\gamma,i}\mid\gamma\in\Gamma,1\leq i\leq k_{\gamma}\}$ be
a nilpotent Jordan normal base of $_{R}M$\ with respect to the nilpotent
endomorphism $\varphi\in$ Hom$_{R}(M,M)$. We keep the notations of the
previous section and in the rest of the paper we assume that $_{R}M$ is
finitely generated, i.e. that $\Gamma$\ is finite.

\noindent A linear order on $\Gamma$ (say $\Gamma=\{1,2,\ldots,m\}$) allows us
to view an element $\mathbf{f}=(f_{\gamma}(t))_{\gamma\in\Gamma}$ of
$(R[t])^{\Gamma}$ as a $1\times\Gamma$ matrix (a row vector) over $R[t]$. For
a $\Gamma\times\Gamma$ matrix $\mathbf{P}=[p_{\delta,\gamma}(t)]$ in
$M_{\Gamma\times\Gamma}(R[t])$ the matrix product%
\[
\mathbf{fP=}\underset{\delta\in\Gamma}{%
{\displaystyle\sum}
}f_{\delta}(t)\mathbf{p}_{\delta},
\]
of $\mathbf{f}$ and $\mathbf{P}$ is a $1\times\Gamma$ matrix (row vector) in
$(R[t])^{\Gamma}$, where $\mathbf{p}_{\delta}=(p_{\delta,\gamma}%
(t))_{\gamma\in\Gamma}$ is the $\delta$-th row vector of $\mathbf{P}$ and%
\[
(\mathbf{fP})_{\gamma}=\underset{\delta\in\Gamma}{%
{\displaystyle\sum}
}f_{\delta}(t)p_{\delta,\gamma}(t).
\]
Using the homomorphism $\Phi:\underset{\gamma\in\Gamma}{%
{\displaystyle\coprod}
}R[t]\longrightarrow M$ introduced in Section 3, we define the subset%
\[
\mathcal{M}(\Phi)=\{\mathbf{P}\in M_{\Gamma\times\Gamma}(R[t])\mid
\mathbf{fP}\in\ker(\Phi)\text{ for all }\mathbf{f}\in\ker(\Phi)\}.
\]
\noindent If $R$ is a local ring, then we can determine $\mathcal{M}(\Phi)$.
Let $\mathbf{e}_{\delta}\in\ker\Phi$ be the vector with $t^{k_{\delta}}$ in
its $\delta$-coordinate and zeros in all other places. If $\mathbf{P}%
\in\mathcal{M}(\Phi)$, then $\mathbf{e}_{\delta}\mathbf{P}\in\ker\Phi$ implies
that $t^{k_{\delta}}p_{\delta,\gamma}(t)\in J(R)[t]+(t^{k_{\gamma}})$. Thus
for local rings we have%
\[
\mathcal{M}(\Phi)\!=\!\{\mathbf{P}\!\in\!M_{\Gamma\times\Gamma}(R[t])\!\mid
\!\mathbf{P}\!=\![p_{\delta,\gamma}(t)]\text{ and }t^{k_{\delta}}%
p_{\delta,\gamma}(t)\!\in\!J(R)[t]\!+\!(t^{k_{\gamma}})\text{ for all }%
\delta,\gamma\!\in\!\Gamma\}
\]
and $\mathbf{E}_{\delta,\gamma}\in\mathcal{M}(\Phi)$ for all $\delta,\gamma
\in\Gamma$ with $k_{\delta}\geq k_{\gamma}$ (where $\mathbf{E}_{\delta,\gamma
}$ denotes the $\Gamma\times\Gamma$\ standard matrix unit over $R[t]$\ with
$1$ in the $(\delta,\gamma)$ entry and zeros in the other entries).

\bigskip

\noindent\textbf{4.1.Lemma.}\textit{ }$\mathcal{M}(\Phi)$\textit{\ is a
}$Z(R)$\textit{-subalgebra of }$M_{\Gamma\times\Gamma}(R[t])$\textit{. For
}$\mathbf{P}\in\mathcal{M}(\Phi)$\textit{ and }$\mathbf{f}=(f_{\gamma
}(t))_{\gamma\in\Gamma}$\textit{ in }$(R[t])^{\Gamma}$\textit{ the formula}%
\[
\psi_{\mathbf{P}}(\Phi(\mathbf{f}))=\Phi(\mathbf{fP})
\]
\textit{properly defines an }$R$\textit{-endomorphism }$\psi_{\mathbf{P}%
}:M\longrightarrow M$\textit{ of }$_{R}M$\textit{ such that }$\psi
_{\mathbf{P}}\circ\varphi=\varphi\circ\psi_{\mathbf{P}}$\textit{. The
assignment }$\mathbf{P}\longmapsto\psi_{\mathbf{P}}$\textit{ is a}%
\[
\mathcal{M}(\Phi)^{\text{op}}\longrightarrow C_{\varphi}%
\]
\textit{homomorphism of }$Z(R)$\textit{-algebras.}

\bigskip

\noindent\textbf{Proof.} For $\mathbf{P},\mathbf{Q}\in\mathcal{M}(\Phi)$,
$\mathbf{f}\in\ker(\Phi)$ and $c\in Z(R)$ we have%
\[
\Phi(\mathbf{f}(c\mathbf{P}))=\Phi(c(\mathbf{fP}))=c\Phi(\mathbf{fP})=0,
\]%
\[
\Phi(\mathbf{f}(\mathbf{P}+\mathbf{Q}))=\Phi(\mathbf{fP}+\mathbf{fQ}%
)=\Phi(\mathbf{fP})+\Phi(\mathbf{fQ})=0
\]
and $\mathbf{fP}\in\ker(\Phi)$ implies that%
\[
\Phi(\mathbf{f}(\mathbf{PQ}))=\Phi((\mathbf{fP})\mathbf{Q})=0,
\]
whence $\mathbf{f}(\mathbf{PQ})\in\ker(\Phi)$ follows. Thus $c\mathbf{P}%
,\mathbf{P}+\mathbf{Q},\mathbf{PQ}\in\mathcal{M}(\Phi)$.

\noindent Let $\mathbf{g}\in(R[t])^{\Gamma}$. If $\Phi(\mathbf{f}%
)=\Phi(\mathbf{g})$, then $\mathbf{f}-\mathbf{g}\in\ker(\Phi)$ implies that
$(\mathbf{f}-\mathbf{g})\mathbf{P}\in\ker(\Phi)$, whence $\Phi(\mathbf{fP}%
)=\Phi(\mathbf{gP})$ follows. Since $\Phi$ is surjective, it follows that
$\psi_{\mathbf{P}}$ is well defined. It is straightforward to check that%
\[
\psi_{\mathbf{P}}(\Phi(\mathbf{f})+\Phi(\mathbf{g}))=\psi_{\mathbf{P}}%
(\Phi(\mathbf{f}))+\psi_{\mathbf{P}}(\Phi(\mathbf{g}))\text{ and }%
\psi_{\mathbf{P}}(r\Phi(\mathbf{f}))=r\psi_{\mathbf{P}}(\Phi(\mathbf{f}))
\]
for all $\mathbf{f},\mathbf{g}\in(R[t])^{\Gamma}$ and $r\in R$. Thus
$\psi_{\mathbf{P}}$ is an $R$-endomorphism. In view of%
\[
\psi_{\mathbf{P}}(\varphi(\Phi(\mathbf{f}))=\psi_{\mathbf{P}}(t\ast
\Phi(\mathbf{f}))=\psi_{\mathbf{P}}(\Phi(t\mathbf{f}))=\Phi((t\mathbf{f}%
)\mathbf{P})=
\]%
\[
=\Phi(t(\mathbf{fP}))=t\ast\Phi(\mathbf{fP})=t\ast\psi_{\mathbf{P}}%
(\Phi(\mathbf{f}))=\varphi(\psi_{\mathbf{P}}(\Phi(\mathbf{f}))),
\]
the surjectivity of $\Phi$\ gives that $\psi_{\mathbf{P}}\circ\varphi
=\varphi\circ\psi_{\mathbf{P}}$. Clearly,%
\[
\psi_{c\mathbf{P}}=c\psi_{\mathbf{P}}\text{ , }\psi_{\mathbf{P}+\mathbf{Q}%
}=\psi_{\mathbf{P}}+\psi_{\mathbf{Q}}\text{ and }\psi_{\mathbf{PQ}}%
=\psi_{\mathbf{Q}}\circ\psi_{\mathbf{P}}%
\]
ensure that $\mathbf{P}\longmapsto\psi_{\mathbf{P}}$ is a homomorphism of
$Z(R)$-algebras. We deal only with the last identity:%
\[
\psi_{\mathbf{PQ}}(\Phi(\mathbf{f}))=\Phi(\mathbf{f}(\mathbf{PQ}%
))=\Phi((\mathbf{fP})\mathbf{Q})=\psi_{\mathbf{Q}}(\Phi(\mathbf{fP}%
))=\psi_{\mathbf{Q}}(\psi_{\mathbf{P}}(\Phi(\mathbf{f}))
\]
proves our claim.$\square$

\bigskip

\noindent\textbf{4.2.Lemma.}\textit{ The }$R$\textit{-submodule}%
\[
\mathcal{V}(k_{\gamma},\mathbf{\!}\gamma\in\Gamma)\mathbf{\!}=\mathbf{\!}%
\{\mathbf{P\!}\in\mathbf{\!}M_{\Gamma\times\Gamma}(R[t])\mathbf{\!}%
\mid\mathbf{\!P\!}=\mathbf{\!}[p_{\delta,\gamma}(t)]\text{\thinspace
and\thinspace}\deg(p_{\delta,\gamma}(t))\mathbf{\!}\leq\mathbf{\!}k_{\gamma
}-\mathbf{\!}1\text{\thinspace for\thinspace all\thinspace}\delta
,\mathbf{\!}\gamma\mathbf{\!}\in\mathbf{\!}\Gamma\}
\]
\textit{of }$M_{\Gamma\times\Gamma}(R[t])$\textit{ is }$R$\textit{-generated
by the set }$\{t^{i}\mathbf{E}_{\delta,\gamma}\mid\delta,\gamma\in\Gamma$ and
$0\leq i\leq k_{\gamma}-1\}$\textit{\ of matrices:}%
\[
\mathcal{V}(k_{\gamma},\gamma\in\Gamma)=\underset{\delta,\gamma\in\Gamma,0\leq
i\leq k_{\gamma}-1}{\sum}Rt^{i}\mathbf{E}_{\delta,\gamma}.
\]

\bigskip

\noindent\textbf{Proof.} If $\mathbf{P}=[p_{\delta,\gamma}(t)]$ and
$\mathbf{P}\in\mathcal{V}(k_{\gamma},\gamma\in\Gamma)$, then%
\[
\mathbf{P}=\underset{\delta,\gamma\in\Gamma}{\sum}p_{\delta,\gamma
}(t)\mathbf{E}_{\delta,\gamma}=\underset{\delta,\gamma\in\Gamma,0\leq i\leq
k_{\gamma}-1}{\sum}b_{\delta,\gamma,i}t^{i}\mathbf{E}_{\delta,\gamma},
\]
where $\deg(p_{\delta,\gamma}(t))\leq k_{\gamma}-1$ and%
\[
p_{\delta,\gamma}(t)=b_{\delta,\gamma,1}+b_{\delta,\gamma,2}t+\cdots
+b_{\delta,\gamma,k_{\gamma}}t^{k_{\gamma}-1}.\square
\]

\bigskip

\noindent\textbf{4.3.Lemma.}\textit{ If }$\psi\circ\varphi=\varphi\circ\psi
$\textit{ holds for an }$R$\textit{-endomorphism }$\psi:M\longrightarrow
M$\textit{ of }$_{R}M$\textit{, then there exists a }$\Gamma\times\Gamma
$\textit{ matrix }$\mathbf{P}\in\mathcal{M}(\Phi)\cap\mathcal{V}(k_{\gamma
},\gamma\in\Gamma)$\textit{ such that}%
\[
\psi(\Phi(\mathbf{f}))=\Phi(\mathbf{fP})
\]
\textit{for all }$\mathbf{f}=(f_{\gamma}(t))_{\gamma\in\Gamma}$\textit{ in
}$(R[t])^{\Gamma}$\textit{.}

\bigskip

\noindent\textbf{Proof.} Since $\Phi:(R[t])^{\Gamma}\longrightarrow M$ is
surjective, for each $\delta\in\Gamma$\ we can find an element $\mathbf{p}%
_{\delta}=(p_{\delta,\gamma}(t))_{\gamma\in\Gamma}$ in $(R[t])^{\Gamma}$ such
that $\Phi(\mathbf{p}_{\delta})=\psi(x_{\delta,1})$. We have $\Phi
(\mathbf{p}_{\delta}^{(c)})=\psi(x_{\delta,1})$, where $\mathbf{p}_{\delta
}^{(c)}=(p_{\delta,\gamma}^{(k_{\gamma})}(t))_{\gamma\in\Gamma}$ is the cut of
$\mathbf{p}_{\delta}$\ with respect to the given nilpotent Jordan normal base.
The $\Gamma\times\Gamma$ matrix $\mathbf{P}=[p_{\delta,\gamma}^{(k_{\gamma}%
)}(t)]$ consisting of the row vectors $\mathbf{p}_{\delta}^{(c)}$, $\delta
\in\Gamma$, is in $\mathcal{V}(k_{\gamma},\gamma\in\Gamma)$ and%
\[
\psi(\Phi(\mathbf{f}))=\underset{\delta\in\Gamma}{%
{\displaystyle\sum}
}\psi(f_{\delta}(t)\ast x_{\delta,1})=\underset{\delta\in\Gamma}{%
{\displaystyle\sum}
}f_{\delta}(t)\ast\psi(x_{\delta,1})=\underset{\delta\in\Gamma}{%
{\displaystyle\sum}
}f_{\delta}(t)\ast\Phi(\mathbf{p}_{\delta}^{(c)})=
\]%
\[
=\underset{\delta\in\Gamma}{%
{\displaystyle\sum}
}\Phi(f_{\delta}(t)\mathbf{p}_{\delta}^{(c)})=\Phi(\underset{\delta\in\Gamma}{%
{\displaystyle\sum}
}f_{\delta}(t)\mathbf{p}_{\delta}^{(c)})=\Phi(\mathbf{fP})
\]
for all $\mathbf{f}\in(R[t])^{\Gamma}$. Since $\mathbf{f}\in\ker(\Phi)$
implies that $\Phi(\mathbf{fP})=\psi(\Phi(\mathbf{f}))=0$, we obtain that
$\mathbf{P}\in\mathcal{M}(\Phi)$.$\square$

\bigskip

\noindent\textbf{4.4.Theorem.}\textit{ Let }$\varphi:M\longrightarrow
M$\textit{ be a nilpotent }$R$\textit{-endomorphism of the finitely generated
semisimple left }$R$\textit{-module }$_{R}M$\textit{. Then the centralizer
}$C_{\varphi}$\textit{\ (as a }$Z(R)$\textit{-subalgebra of }Hom$_{R}%
(M,M)$\textit{) is the homomorphic image of the opposite of some }%
$Z(R)$\textit{-subalgebra of the matrix algebra }$M_{m\times m}(R[t])$%
\textit{, where }$m=\dim_{R}(\ker(\varphi))$\textit{.}

\bigskip

\noindent\textbf{Proof.} Lemma 4.1 ensures that $\mathbf{P}\longmapsto
\psi_{\mathbf{P}}$ is a%
\[
\mathcal{M}(\Phi)^{\text{op}}\longrightarrow C_{\varphi}%
\]
homomorphism of $Z(R)$-algebras. Since our assignment is surjective by Lemma
4.3, we obtain that $C_{\varphi}$ is the homomorphic image of $\mathcal{M}%
(\Phi)^{\text{op}}$, where $\mathcal{M}(\Phi)$ is a $Z(R)$-subalgebra of
$M_{\Gamma\times\Gamma}(R[t])$. To conclude the proof it is enough to note
that $m=\dim_{R}(\ker(\varphi))=\left\vert \Gamma\right\vert $.$\square$

\bigskip

\noindent\textbf{4.5.Corollary.}\textit{ Let }$\varphi:M\longrightarrow
M$\textit{ be a nilpotent }$R$\textit{-endomorphism of the finitely generated
semisimple left }$R$\textit{-module }$_{R}M$\textit{. Then }$C_{\varphi}%
$\textit{ satisfies all of the polynomial identities (with coefficients in
}$Z(R)$\textit{) of }$M_{m\times m}^{\text{op}}(R[t])$\textit{. If }%
$R$\textit{\ is commutative, then }$C_{\varphi}$\textit{ satisfies the
standard identity }$S_{2m}=0$\textit{ of degree }$2m$\textit{ by the
Amitsur-Levitzki theorem.}

\bigskip

\noindent\textbf{4.6.Theorem.}\textit{ Let }$\varphi:M\longrightarrow
M$\textit{ be a nilpotent }$R$\textit{-endomorphism of the finitely generated
semisimple left }$R$\textit{-module }$_{R}M$\textit{. Then the centralizer
}$C_{\varphi}$\textit{\ (as a }$Z(R)$\textit{-submodule of }Hom$_{R}%
(M,M)$\textit{) is the homomorphic image of some }$Z(R)$\textit{-submodule of
a certain }$md$\textit{-generated\ }$R$\textit{-submodule of }$M_{m\times
m}(R[t])$\textit{, where }$d=\dim_{R}(M)$\textit{ and }$m=\dim_{R}%
(\ker(\varphi))$\textit{.}

\bigskip

\noindent\textbf{Proof.} Lemma 4.1 and 4.3 ensure that $\mathbf{P}%
\longmapsto\psi_{\mathbf{P}}$ is a surjective%
\[
\mathcal{M}(\Phi)\cap\mathcal{V}(k_{\gamma},\gamma\in\Gamma)\longrightarrow
C_{\varphi}%
\]
homomorphism of $Z(R)$-modules, where $\mathcal{M}(\Phi)\cap\mathcal{V}%
(k_{\gamma},\gamma\in\Gamma)$ is a $Z(R)$-submodule of the left $R$-submodule
$\mathcal{V}(k_{\gamma},\gamma\in\Gamma)$ of $M_{\Gamma\times\Gamma}(R[t])$.
Using Lemma 4.2, we obtain that $\mathcal{V}(k_{\gamma},\gamma\in\Gamma)$ is
$R$-generated by the set $\{t^{i}\mathbf{E}_{\delta,\gamma}\mid\delta
,\gamma\in\Gamma$ and $0\leq i\leq k_{\gamma}-1\}$ having%
\[
\left\vert \Gamma\right\vert \cdot\underset{\gamma\in\Gamma}{\sum}k_{\gamma
}=md
\]
elements.$\square$

\bigskip

\noindent\textbf{4.7.Corollary.}\textit{ Let }$\varphi:M\longrightarrow
M$\textit{ be a nilpotent }$D$\textit{-linear map of the finite dimensional
left vector space }$_{D}M$\textit{\ over a division ring }$D$\textit{. If }%
$D$\textit{ is finite dimensional over its central subfield }$K=Z(D)$\textit{,
then }$C_{\varphi}$\textit{\ is a\ }$K$\textit{-subspace of }Hom$_{D}%
(M,M)$\textit{\ and}%
\[
\dim_{K}(C_{\varphi})\leq\dim_{K}(D)\cdot\dim_{D}(\ker(\varphi))\cdot\dim
_{D}(M).
\]
\noindent\textbf{Remark.} If $D=K$ is an algebraically closed field and all
Jordan blocks of the nilpotent linear map $\varphi\in$Hom$_{K}(M,M)$\ (of the
finite dimensional $K$-vector space $M$) are of the same size $s\times s$,
then we can use the formula for $\dim(C_{A_{i}})$\ in Section 1 to see that
the upper bound in the above Corollary 4.7 is sharp:%
\[
\dim_{K}(C_{\varphi})=m^{2}s=md=\dim_{K}(\ker(\varphi))\cdot\dim_{K}(M).
\]

\bigskip

\noindent\textbf{4.8.Theorem.}\textit{ If }$_{R}M$\textit{ is semisimple and
}$\varphi:M\longrightarrow M$\textit{ is an indecomposable nilpotent element
of the ring }Hom$_{R}(M,M)$\textit{, then the following are equivalent.}

\begin{enumerate}
\item $\psi\in C_{\varphi}$\textit{.}

\item \textit{We can find an }$R$\textit{-generating set }$\{y_{j}\in
M\mid1\leq j\leq d\}$\textit{ of }$_{R}M$\textit{ and elements }$a_{1}%
,a_{2},\ldots,a_{n}$\textit{ in }$R$\textit{ such that}%
\[
a_{1}y_{j}+a_{2}\varphi(y_{j})+\cdots+a_{n}\varphi^{n-1}(y_{j})=\psi(y_{j})
\]
\textit{and}%
\[
a_{1}\varphi(y_{j})+a_{2}\varphi(\varphi(y_{j}))+\cdots+a_{n}\varphi
^{n-1}(\varphi(y_{j}))=\psi(\varphi(y_{j}))
\]
\textit{for all }$1\leq j\leq d$\textit{.}
\end{enumerate}

\bigskip

\noindent\textbf{Proof.}

\noindent(1)$\Longrightarrow$(2): Obviously, if $\psi\in C_{\varphi}$ then the
first identity implies the second one. Proposition 2.3 ensures the existence
of a nilpotent Jordan normal base $\{x_{i}\mid1\leq i\leq n\}$\ of $_{R}%
M$\ with respect to $\varphi$ consisting of one block. Clearly, $\underset
{1\leq i\leq n}{\oplus}Rx_{i}=M$ implies that%
\[
\psi(x_{1})=a_{1}x_{1}+a_{2}x_{2}+\cdots+a_{n}x_{n}=a_{1}x_{1}+a_{2}%
\varphi(x_{1})+\cdots+a_{n}\varphi^{n-1}(x_{1})
\]
for some $a_{1},a_{2},\ldots,a_{n}\in R$. Thus%
\[
\psi(x_{i})=\psi(\varphi^{i-1}(x_{1}))=\varphi^{i-1}(\psi(x_{1}))=\varphi
^{i-1}(a_{1}x_{1}+a_{2}\varphi(x_{1})+\cdots+a_{n}\varphi^{n-1}(x_{1}))=
\]%
\[
=\!a_{1}\varphi^{i-1}(x_{1})\!+\!a_{2}\varphi(\varphi^{i-1}(x_{1}%
))\!+\cdots+\!a_{n}\varphi^{n-1}(\varphi^{i-1}(x_{1}))\!=\!a_{1}%
x_{i}\!+\!a_{2}\varphi(x_{i})\!+\cdots+\!a_{n}\varphi^{n-1}(x_{i})
\]
for all $1\leq i\leq n$.

\noindent(2)$\Longrightarrow$(1): Since we have%
\[
\varphi(\psi(y_{j}))=\varphi(a_{1}y_{j}+a_{2}\varphi(y_{j})+\cdots
+a_{n}\varphi^{n-1}(y_{j}))=
\]%
\[
=a_{1}\varphi(y_{j})+a_{2}\varphi(\varphi(y_{j}))+\cdots+a_{n}\varphi
^{n-1}(\varphi(y_{j}))=\psi(\varphi(y_{j}))
\]
for all $1\leq j\leq d$, the implication is proved.$\square$

\bigskip

\noindent\textbf{4.9.Corollary.}\textit{ If }$R$\textit{ is commutative,
}$_{R}M$\textit{ is semisimple and }$\varphi:M\longrightarrow M$\textit{ is an
indecomposable nilpotent element of the ring }Hom$_{R}(M,M)$\textit{, then the
following are equivalent.}

\begin{enumerate}
\item $\psi\in C_{\varphi}$\textit{.}

\item \textit{We can find elements }$a_{1},a_{2},\ldots,a_{n}$\textit{ in }%
$R$\textit{ such that}%
\[
a_{1}u+a_{2}\varphi(u)+\cdots+a_{n}\varphi^{n-1}(u)=\psi(u)
\]
\textit{for all }$u\in M$\textit{. In other words }$\psi$\textit{\ is a
polynomial of }$\varphi$\textit{.}
\end{enumerate}

\bigskip

\noindent\textbf{Proof.} It suffices to prove that if $\underset{1\leq j\leq
d}{\sum}Ry_{j}=M$ and%
\[
a_{1}y_{j}+a_{2}\varphi(y_{j})+\cdots+a_{n}\varphi^{n-1}(y_{j})=\psi(y_{j})
\]
holds for all $1\leq j\leq d$, then we have%
\[
a_{1}u+a_{2}\varphi(u)+\cdots+a_{n}\varphi^{n-1}(u)=\psi(u)
\]
for all $u\in M$. Since $u=b_{1}y_{1}+b_{2}y_{2}+\cdots+b_{d}y_{d}$ for some
$b_{1},b_{2},\ldots,b_{d}\in R$ and $b_{j}a_{i}=a_{i}b_{j}$, we obtain that%
\[
\psi(u)=\underset{1\leq j\leq d}{\sum}b_{j}\psi(y_{j})=\underset{1\leq j\leq
d}{\sum}b_{j}(a_{1}y_{j}+a_{2}\varphi(y_{j})+\cdots+a_{n}\varphi^{n-1}%
(y_{j}))=
\]%
\[
=a_{1}\left(  \underset{1\leq j\leq d}{\sum}b_{j}y_{j}\right)  +a_{2}%
\varphi\left(  \underset{1\leq j\leq d}{\sum}b_{j}y_{j}\right)  +\cdots
+a_{n}\varphi^{n-1}\left(  \underset{1\leq j\leq d}{\sum}b_{j}y_{j}\right)  =
\]%
\[
=a_{1}u+a_{2}\varphi(u)+\cdots+a_{n}\varphi^{n-1}(u).\square
\]

\bigskip

\noindent\textbf{4.10.Theorem.}\textit{ Let }$R$\textit{\ be a local ring. If
}$\varphi:M\longrightarrow M$\textit{ is a nilpotent }$R$\textit{-endomorphism
of the finitely generated semisimple left }$R$\textit{-module }$_{R}M$\textit{
and}

\noindent$\sigma\in$ Hom$_{R}(M,M)$\textit{ is arbitrary, then the following
are equivalent.}

\begin{enumerate}
\item $C_{\varphi}\subseteq C_{\sigma}$\textit{.}

\item \textit{We can find an }$R$\textit{-generating set }$\{y_{j}\in
M\mid1\leq j\leq d\}$\textit{ of }$_{R}M$\textit{ and elements }$a_{1}%
,a_{2},\ldots,a_{n}$\textit{ in }$R$\textit{ such that}%
\[
a_{1}\psi(y_{j})+a_{2}\varphi(\psi(y_{j}))+\cdots+a_{n}\varphi^{n-1}%
(\psi(y_{j}))=\sigma(\psi(y_{j}))
\]
\textit{for all }$1\leq j\leq d$\textit{ and all }$\psi\in C_{\varphi}%
$\textit{.}
\end{enumerate}

\bigskip

\noindent\textbf{Proof.}

\noindent(1)$\Longrightarrow$(2): Obviously, if $C_{\varphi}\subseteq
C_{\sigma}$ then%
\[
a_{1}y_{j}+a_{2}\varphi(y_{j})+\cdots+a_{n}\varphi^{n-1}(y_{j})=\sigma(y_{j})
\]
implies that%
\[
a_{1}\psi(y_{j})+a_{2}\varphi(\psi(y_{j}))+\cdots+a_{n}\varphi^{n-1}%
(\psi(y_{j}))=\sigma(\psi(y_{j}))
\]
for all $\psi\in C_{\varphi}$. Theorem 2.1 ensures the existence of a
nilpotent Jordan normal base $\{x_{\gamma,i}\mid\gamma\in\Gamma,1\leq i\leq
k_{\gamma}\}$ of $_{R}M$\ with respect to $\varphi$. Consider the natural
projection $\varepsilon_{\delta}:M\longrightarrow N_{\delta}^{\prime}$
corresponding to the direct sum $M=N_{\delta}^{\prime}\oplus N_{\delta
}^{\prime\prime}$ (see the proof of 2.3), where
\[
N_{\delta}^{\prime}=\underset{1\leq i\leq k_{\delta}}{\oplus}Rx_{\delta
,i}\text{ and }N_{\delta}^{\prime\prime}=\underset{\gamma\in\Gamma
\setminus\{\delta\},1\leq i\leq k_{\gamma}}{\oplus}Rx_{\gamma,i}\text{ }.
\]
Then $\varepsilon_{\delta}\in C_{\varphi}$, whence $\varepsilon_{\delta}\in
C_{\sigma}$ follows for all $\delta\in\Gamma$. Thus im$(\varepsilon_{\delta
})=N_{\delta}^{\prime}$ and%
\[
\sigma:\text{im}(\varepsilon_{\delta})\longrightarrow\text{im}(\varepsilon
_{\delta})
\]
implies that%
\[
\sigma(x_{\delta,1})=\underset{1\leq i\leq k_{\delta}}{%
{\displaystyle\sum}
}a_{\delta,i}x_{\delta,i}=h_{\delta}(t)\ast x_{\delta,1}%
\]
for some $h_{\delta}(t)=a_{\delta,1}+a_{\delta,2}t+\cdots+a_{\delta,k_{\delta
}}t^{k_{\delta}-1}$ in $R[t]$. Since $\varphi\in C_{\sigma}$ implies that
$\sigma\in C_{\varphi}$, it follows that%
\[
\sigma(\Phi(\mathbf{f}))=\underset{\gamma\in\Gamma}{%
{\displaystyle\sum}
}\sigma(f_{\gamma}(t)\ast x_{\gamma,1})=\underset{\gamma\in\Gamma}{%
{\displaystyle\sum}
}f_{\gamma}(t)\ast\sigma(x_{\gamma,1})=\underset{\gamma\in\Gamma}{%
{\displaystyle\sum}
}f_{\gamma}(t)\ast(h_{\gamma}(t)\ast x_{\gamma,1})=
\]%
\[
=\underset{\gamma\in\Gamma}{%
{\displaystyle\sum}
}(f_{\gamma}(t)h_{\gamma}(t))\ast x_{\gamma,1}=\Phi(\mathbf{fH}),
\]
where $\mathbf{f}\in(R[t])^{\Gamma}$ and $\mathbf{H}=\underset{\gamma\in
\Gamma}{\sum}h_{\gamma}(t)\mathbf{E}_{\gamma,\gamma}$ is a $\Gamma\times
\Gamma$ diagonal matrix in $\mathcal{M}(\Phi)$ (note that $\mathbf{H\in
}\mathcal{M}(\Phi)$ is a consequence of $\sigma(\Phi(\mathbf{f}))=\Phi
(\mathbf{fH})$). In view of Lemma 4.1 and 4.3, the containment $C_{\varphi
}\subseteq C_{\sigma}$ is equivalent to the condition that $\sigma\circ
\psi_{\mathbf{P}}=\psi_{\mathbf{P}}\circ\sigma$ for all $\mathbf{P}%
\in\mathcal{M}(\Phi)$. Consequently, we obtain that $C_{\varphi}\subseteq
C_{\sigma}$ is equivalent to the following:%
\[
\Phi(\mathbf{fPH})=\sigma(\Phi(\mathbf{fP})))=\sigma(\psi_{\mathbf{P}}%
(\Phi(\mathbf{f})))=\psi_{\mathbf{P}}(\sigma(\Phi(\mathbf{f})))=\psi
_{\mathbf{P}}(\Phi(\mathbf{fH}))=\Phi(\mathbf{fHP})
\]
for all $\mathbf{f}\in(R[t])^{\Gamma}$ and $\mathbf{P}\in\mathcal{M}(\Phi)$.
Thus $C_{\varphi}\subseteq C_{\sigma}$ implies that $\Phi(\mathbf{f}%
(\mathbf{PH}-\mathbf{HP}))=0$, i.e. that $\mathbf{f}(\mathbf{PH}%
-\mathbf{HP})\in\ker(\Phi)$ for all $\mathbf{f}$ and $\mathbf{P}$.

\noindent Now we use that $R$\ is a local ring. If $\alpha\in\Gamma$ is an
index such that%
\[
k_{\alpha}=\max\{k_{\gamma}\mid\gamma\in\Gamma\}=n,
\]
then $\mathbf{E}_{\alpha,\delta}\in\mathcal{M}(\Phi)$ for all $\delta\in
\Gamma$ (see the argument preceding Lemma 4.1). Take $\mathbf{e}%
=(1)_{\gamma\in\Gamma}$ and $\mathbf{P}=\mathbf{E}_{\alpha,\delta}$, then the
$\delta$-coordinate of%
\[
\mathbf{e}(\mathbf{E}_{\alpha,\delta}\mathbf{H}-\mathbf{HE}_{\alpha,\delta
})=(h_{\delta}(t)-h_{\alpha}(t))\mathbf{eE}_{\alpha,\delta}%
\]
is $h_{\delta}(t)-h_{\alpha}(t)$. Since%
\[
(h_{\delta}(t)-h_{\alpha}(t))\mathbf{eE}_{\alpha,\delta}\in\ker(\Phi
)=\underset{\gamma\in\Gamma}{%
{\displaystyle\coprod}
}J(R)[t]+(t^{k_{\gamma}}),
\]
we obtain that $h_{\delta}(t)-h_{\alpha}(t)\in J(R)[t]+(t^{k_{\delta}})$. Thus%
\[
\sigma(x_{\delta,1})=h_{\delta}(t)\ast x_{\delta,1}=h_{\alpha}(t)\ast
x_{\delta,1}%
\]
for all $\delta\in\Gamma$. It follows that%
\[
\sigma(x_{\gamma,i})=\sigma(\varphi^{i-1}(x_{\gamma,1}))=\varphi^{i-1}%
(\sigma(x_{\gamma,1}))=\varphi^{i-1}(h_{\alpha}(t)\ast x_{\gamma,1})=
\]%
\[
=h_{\alpha}(t)\ast\varphi^{i-1}(x_{\gamma,1})=h_{\alpha}(t)\ast x_{\gamma
,i}=a_{1}x_{\gamma,i}+a_{2}\varphi(x_{\gamma,i})+\cdots+a_{n}\varphi
^{n-1}(x_{\gamma,i}),
\]
where $h_{\alpha}(t)=a_{1}+a_{2}t+\cdots+a_{n}t^{n-1}$.

\noindent(2)$\Longrightarrow$(1): ($R$ is an arbitrary ring) Since $1_{M}\in
C_{\varphi}$, we obtain that%
\[
a_{1}y_{j}+a_{2}\varphi(y_{j})+\cdots+a_{n}\varphi^{n-1}(y_{j})=\sigma(y_{j})
\]
for all $1\leq j\leq d$. If $\psi\in C_{\varphi}$, then%
\[
\psi(\sigma(y_{j}))=\psi(a_{1}y_{j}+a_{2}\varphi(y_{j})+\cdots+a_{n}%
\varphi^{n-1}(y_{j}))=
\]%
\[
=a_{1}\psi(y_{j})+a_{2}\varphi(\psi(y_{j}))+\cdots+a_{n}\varphi^{n-1}%
(\psi(y_{j}))=\sigma(\psi(y_{j}))
\]
for all $1\leq j\leq d$, whence $\psi\circ\sigma=\sigma\circ\psi$ follows.
Thus $C_{\varphi}\subseteq C_{\sigma}$.$\square$

\bigskip

\noindent5. THE\ CENTRALIZER\ OF\ AN ARBITRARY LINEAR MAP

\bigskip

\noindent Let $r\geq1$ be an integer such that%
\[
\ker(\varphi^{r})\oplus\text{im}(\varphi^{r})=M
\]
for the $R$-endomorphism $\varphi\in$Hom$_{R}(M,M)$ of the left $R$-module
$_{R}M$. Since%
\[
\varphi:\ker(\varphi^{r})\longrightarrow\ker(\varphi^{r})\text{ and }%
\varphi:\text{im}(\varphi^{r})\longrightarrow\text{im}(\varphi^{r}),
\]
we can take the restricted $R$-endomorphisms%
\[
\varphi_{1}=\varphi\upharpoonright\ker(\varphi^{r})\text{ and }\varphi
_{2}=\varphi\upharpoonright\text{im}(\varphi^{r}).
\]
In the next statement we consider the following centralizers%
\[
C_{\varphi}\subseteq\text{Hom}_{R}(M,M),C_{\varphi_{1}}\subseteq\text{Hom}%
_{R}(\ker(\varphi^{r}),\ker(\varphi^{r})),C_{\varphi_{2}}\subseteq
\text{Hom}_{R}(\text{im}(\varphi^{r}),\text{im}(\varphi^{r})).
\]

\bigskip

\noindent\textbf{5.1.Lemma.}\textit{ If }$\ker(\varphi^{r})$ $\oplus$
im$(\varphi^{r})=M$\textit{ holds for }$\varphi\in$Hom$_{R}(M,M)$\textit{,
then we have an isomorphism}%
\[
C_{\varphi}\cong C_{\varphi_{1}}\times C_{\varphi_{2}}%
\]
\textit{of }$Z(R)$\textit{-algebras.}

\bigskip

\noindent\textbf{Proof.} Consider the natural projections%
\[
\varepsilon_{1}:M\longrightarrow\ker(\varphi^{r})\text{ and }\varepsilon
_{2}:M\longrightarrow\text{im}(\varphi^{r})
\]
and the natural injections%
\[
\tau_{1}:\ker(\varphi^{r})\longrightarrow M\text{ and }\tau_{2}:\text{im}%
(\varphi^{r})\longrightarrow M.
\]
For $\sigma\in C_{\varphi}$ define an assignment by%
\[
\sigma\longmapsto(\varepsilon_{1}\sigma\tau_{1},\varepsilon_{2}\sigma\tau_{2})
\]
and for a pair $(\sigma_{1},\sigma_{2})\in C_{\varphi_{1}}\times
C_{\varphi_{2}}$ define an assignment by%
\[
(\sigma_{1},\sigma_{2})\longmapsto\tau_{1}\sigma_{1}\varepsilon_{1}+\tau
_{2}\sigma_{2}\varepsilon_{2}.
\]
Since $\sigma\in C_{\varphi}$ ensures that%
\[
\sigma(\ker(\varphi^{r}))\subseteq\ker(\varphi^{r})\text{ and }\sigma
(\text{im}(\varphi^{r}))\subseteq\text{im}(\varphi^{r}),
\]
it is straightforward to see that the above definitions provide%
\[
C_{\varphi}\longrightarrow C_{\varphi_{1}}\times C_{\varphi_{2}}\text{ and
}C_{\varphi_{1}}\times C_{\varphi_{2}}\longrightarrow C_{\varphi}%
\]
homomorphisms of $Z(R)$-algebras which are mutual inverses of each
other.$\square$

\bigskip

\noindent\textbf{5.2.Theorem.}\textit{ Let }$K$\textit{ be an algebraically
closed field and }$\varphi:V\longrightarrow V$\textit{ be a }$K$%
\textit{-linear map of the finite dimensional vector space }$V$\textit{. If}%
\[
\{\lambda_{1},\lambda_{2},\ldots,\lambda_{p}\}\subseteq K
\]
\textit{is the (}$p$\textit{-element) set of all eigenvalues of }$\varphi
$\textit{ and}%
\[
m_{i}=\dim(\ker(\varphi-\lambda_{i}1_{V})),
\]
\textit{then for each }$1\leq i\leq p$\textit{ there exists a }$K$%
\textit{-subalgebra }$\mathcal{M}_{i}$ \textit{of the full }$m_{i}\times
m_{i}$\textit{ matrix algebra }$M_{m_{i}\times m_{i}}(K[t])$\textit{ such that
the centralizer }$C_{\varphi}$\textit{ is the homomorphic image of the direct
product }$K$\textit{-algebra}%
\[
\mathcal{M}_{1}\times\mathcal{M}_{2}\times\cdots\times\mathcal{M}_{p}.
\]
\noindent\textbf{Proof.} We apply induction on the dimension of $V$.

\noindent If $\dim V=1$, then $V=Kx$ and Hom$_{K}(V,V)\cong K$ (as
$K$-algebras) is commutative. Now $\varphi=\lambda1_{V}$ and $\dim
(\ker(\varphi-\lambda1_{V}))=1$, where $\lambda$ is the (only) eigenvalue of
$\varphi$. Thus $C_{\varphi}=$Hom$_{K}(V,V)\cong M_{1\times1}(K)$ and
$M_{1\times1}(K)$ is a $K$-subalgebra of $M_{1\times1}(K[t])$.

\noindent Let $n\geq2$ be an integer and assume that our theorem holds for all
vector spaces over $K$ of dimension less or equal than $n-1$. Consider the
situation described in the theorem with $\dim V=n$. Since $K$\ is
algebraically closed, we have a (non empty) $p$-element set%
\[
\{\lambda_{1},\lambda_{2},\ldots,\lambda_{p}\}\subseteq K
\]
of eigenvalues of $\varphi$. Let $0\neq u\in V$ be an eigenvector with
$\varphi(u)=\lambda_{1}u$. The Fitting lemma ensures the existence of an
integer $r\geq1$ such that%
\[
\ker(\varphi-\lambda_{1}1_{V})^{r}\oplus\text{im}(\varphi-\lambda_{1}%
1_{V})^{r}=V.
\]
Clearly,%
\[
\varphi-\lambda_{1}1_{V}:\ker(\varphi-\lambda_{1}1_{V})^{r}\longrightarrow
\ker(\varphi-\lambda_{1}1_{V})^{r},
\]
thus the restricted function%
\[
\varphi_{1}=(\varphi-\lambda_{1}1_{V})\upharpoonright\ker(\varphi-\lambda
_{1}1_{V})^{r}%
\]
is a nilpotent $K$-linear map of the $K$-subspace $\ker(\varphi-\lambda
_{1}1_{V})^{r}\leq V$. The application of Theorem 4.4 provides a
$K$-subalgebra $\mathcal{N}_{1}$\ of the matrix algebra $M_{m_{1}\times m_{1}%
}(K[t])$ with $m_{1}=\dim(\ker(\varphi_{1}))$ such that the centralizer%
\[
C_{\varphi_{1}}\subseteq\text{Hom}_{K}(\ker(\varphi-\lambda_{1}1_{V})^{r}%
,\ker(\varphi-\lambda_{1}1_{V})^{r})
\]
is the $K$-homomorphic image of the opposite $K$-algebra $\mathcal{N}%
_{1}^{\text{op}}$. Using%
\[
\ker(\varphi-\lambda_{1}1_{V})\subseteq\ker(\varphi-\lambda_{1}1_{V})^{r},
\]
we obtain that $\ker(\varphi_{1})=\ker(\varphi-\lambda_{1}1_{V})$.

\noindent Since $(\varphi-\lambda_{1}1_{V})(u)=0$ implies that $u\in
\ker(\varphi-\lambda_{1}1_{V})^{r}$, we obtain that im$(\varphi-\lambda
_{1}1_{V})^{r}$ is a proper $K$-subspace of $V$, i.e. that im$(\varphi
-\lambda_{1}1_{V})^{r}\neq V$ and hence $\dim\left(  \text{im}(\varphi
-\lambda_{1}1_{V})^{r}\right)  \leq n-1$.

\noindent In view of%
\[
\varphi-\lambda_{1}1_{V}:\text{im}(\varphi-\lambda_{1}1_{V})^{r}%
\longrightarrow\text{im}(\varphi-\lambda_{1}1_{V})^{r},
\]
the application of the induction hypothesis on the restricted $K$-linear map%
\[
\varphi_{2}=(\varphi-\lambda_{1}1_{V})\upharpoonright\text{im}(\varphi
-\lambda_{1}1_{V})^{r}%
\]
yields the existence of $K$-subalgebras $\mathcal{N}_{j}\leq M_{m_{j}^{\prime
}\times m_{j}^{\prime}}(K[t])$, $2\leq j\leq q$, with%
\[
m_{j}^{\prime}=\dim(\ker(\varphi_{2}-\mu_{j}1_{\text{im}(\varphi-\lambda
_{1}1_{V})^{r}}))
\]
such that the centralizer%
\[
C_{\varphi_{2}}\subseteq\text{Hom}_{K}(\text{im}(\varphi-\lambda_{1}1_{V}%
)^{r},\text{im}(\varphi-\lambda_{1}1_{V})^{r})
\]
is the homomorphic image of the direct product $K$-algebra%
\[
\mathcal{N}_{2}\times\mathcal{N}_{3}\times\cdots\times\mathcal{N}_{q}.
\]
Note that%
\[
\{\mu_{2},\mu_{3},\ldots,\mu_{q}\}\subseteq K
\]
is the set of all eigenvalues of $\varphi_{2}$. Now%
\[
\ker(\varphi_{2}\!-\!\mu_{j}1_{\text{im}(\varphi-\lambda_{1}1_{V})^{r}%
})\!=\!\ker((\varphi\!-\!\lambda_{1}1_{V})\!-\!\mu_{j}1_{V})\!\cap
\text{im}(\varphi\!-\!\lambda_{1}1_{V})^{r}\!\subseteq\!\ker((\varphi
\!-\!(\lambda_{1}\!+\!\mu_{j})1_{V})
\]
implies that%
\[
m_{j}^{\prime}=\dim(\ker(\varphi_{2}-\mu_{j}1_{\text{im}(\varphi-\lambda
_{1}1_{V})^{r}}))\leq\dim(\ker((\varphi-(\lambda_{1}+\mu_{j})1_{V}))=m_{i(j)},
\]
where $\lambda_{1}+\mu_{j}=\lambda_{i(j)}$ for some unique $2\leq i(j)\leq p$.
Indeed, each $\lambda_{1}+\mu_{j}$, $2\leq j\leq q$, is an eigenvalue of
$\varphi$ and $\mu_{j}=0$ would imply that $(\varphi-\lambda_{1}%
1_{V})(v)=\varphi_{2}(v)=0$ for some $0\neq v\in$im$(\varphi-\lambda_{1}%
1_{V})^{r}$ in contradiction with%
\[
\ker(\varphi-\lambda_{1}1_{V})\cap\text{im}(\varphi-\lambda_{1}1_{V}%
)^{r}\subseteq\ker(\varphi-\lambda_{1}1_{V})^{r}\cap\text{im}(\varphi
-\lambda_{1}1_{V})^{r}=\{0\}.
\]
Thus $C_{\varphi_{1}}\times C_{\varphi_{2}}$ is the homomorphic image of the
direct product $K$-algebra%
\[
\mathcal{N}_{1}^{\text{op}}\times\mathcal{N}_{2}\times\cdots\times
\mathcal{N}_{q}%
\]
and $m_{j}^{\prime}\leq m_{i(j)}$\ allows us to view the $K$-subalgebra
$\mathcal{N}_{j}$ of $M_{m_{j}^{\prime}\times m_{j}^{\prime}}(K[t])$ as a
$K$-subalgebra of $M_{m_{i(j)}\times m_{i(j)}}(K[t])$. Using Proposition 5.1,
we deduce that the centralizer $C_{\varphi-\lambda_{1}1_{V}}=C_{\varphi}$ is
the homomorphic image of the following direct product of $K$-algebras%
\[
\mathcal{M}_{1}\times\mathcal{M}_{2}\times\cdots\times\mathcal{M}_{p},
\]
where $\mathcal{M}_{1}=\mathcal{N}_{1}^{\text{op}}$, $\mathcal{M}%
_{i(j)}=\mathcal{N}_{j}$ for $2\leq j\leq q$ and $\mathcal{M}_{i}=\{0\}$ if%
\[
i\in\{1,2,\ldots,p\}\setminus\{i(2),i(3),\ldots,i(q)\}.
\]
Since $K[t]$\ is commutative, the transpose involution ensures that
$M_{m_{1}\times m_{1}}^{\text{op}}(K[t])\cong M_{m_{1}\times m_{1}}(K[t])$ as
$K$-algebras. Thus the opposite $K$-algebra $\mathcal{N}_{1}^{\text{op}}$ also
can be considered as a $K$-subalgebra of $M_{m_{1}\times m_{1}}(K[t])$
consisting of the transposed matrices: $\mathcal{N}_{1}^{\text{op}}%
\cong\mathcal{N}_{1}^{\top}$.$\square$

\bigskip

\noindent\textbf{5.3.Corollary.}\textit{ Let }$K$\textit{ be an algebraically
closed field and }$\varphi:V\longrightarrow V$\textit{ be a }$K$%
\textit{-linear map of the finite dimensional vector space }$V$\textit{. If}%
\[
\{\lambda_{1},\lambda_{2},\ldots,\lambda_{p}\}\subseteq K
\]
\textit{is the (}$p$\textit{-element) set of all eigenvalues of }$\varphi
$\textit{ and}%
\[
m=\max\{\dim(\ker(\varphi-\lambda_{i}1_{V}))\mid1\leq i\leq p\}
\]
\textit{then the centralizer }$C_{\varphi}$\textit{ satisfies the polynomial
identities of the full }$m\times m$\textit{ matrix algebra }$M_{m\times
m}(K[t])$\textit{.}

\bigskip

\noindent\textbf{Proof.} Since $m_{i}=\dim(\ker(\varphi-\lambda_{i}1_{V}))\leq
m$, the $m_{i}\times m_{i}$ matrix algebra $M_{m_{i}\times m_{i}}\!(\!K[t]\!)$
satisfies the polynomial identities of the $m\times m$ matrix algebra
$M_{m\times m}(K[t])$. Thus each $\mathcal{M}_{i}$ satisfies the polynomial
identities of the $m\times m$ matrix algebra $M_{m\times m}(K[t])$, whence we
obtain that the same holds for any homomorphic image of their direct
product.$\square$

\bigskip

\noindent\textbf{Remark.} If $A$ is an $n\times n$ matrix over an
algebraically closed field $K$\ with eigenvalues $\lambda_{1},\lambda
_{2},\ldots,\lambda_{p}$, then $\dim(\ker(A-\lambda_{i}I))$ is the number of
blocks in the Jordan normal form of $A$ containing $\lambda_{i}$\ in the
diagonal. Since Hom$_{K}(V,V)\cong M_{n\times n}(K)$, Corollary 5.3 says that
from a PI point of view the centralizer $C_{A}\subseteq M_{n\times n}(K)$\ of
$A$ behaves like a matrix ring (over the polynomial ring $K[t]$) of size much
smaller than $n$. If the characteristic polynomial of $A$ coincides with the
minimal polynomial, then $m_{i}=\dim(\ker(A-\lambda_{i}I))=1$ for each $1\leq
i\leq p$. Thus $m=\max\{m_{i}\mid1\leq i\leq p\}=1$ and Corollary 5.3 gives
the commutativity of $C_{A}$. For such $A$ we have%
\[
C_{A}=\{f(A)\mid f(t)\in K[t]\},
\]
whence the commutativity of $C_{A}$\ also follows (more details in Section 1).

\bigskip

\noindent REFERENCES

\bigskip

\begin{enumerate}
\item Bergman, G.: \textit{Centralizers in free associative algebras}, Trans.
Amer. Math. Soc. 137 (1969), 327-344.

\item Guralnick, R.: \textit{A note on commuting pairs of matrices}, Linear
and Multilinear Algebra 31 (1992), 71-75.

\item Guralnick, R., Sethuraman, B.A.: \textit{Commuting pairs and triples of
matrices and related varieties}, Linear Algebra Appl. 310 (2000), 139--148.

\item Nelson, G.C., Ton-That, T.:\textit{ Multiplicatively closed bases for
C(A)}, Note di Matematica 26, n. 2 (2006), 81--104.

\item Prasolov, V.V.: \textit{Problems and Theorems in Linear Algebra}, Vol.
134 of Translation of Mathematical Monographs, American Mathematical Society,
Providence, Rhode Island, 1994.

\item Suprunenko, D.A. and Tyshkevich, R.I.:\textit{ Commutative Matrices},
Academic Press, New York and London, 1968.

\item Szigeti, J.: \textit{Linear algebra in lattices and nilpotent
endomorphisms of semisimple modules}, Journal of Algebra 319 (2008), 296--308.
\end{enumerate}

\end{document}